\documentclass[11pt]{article}

\usepackage[a4paper,margin=1in]{geometry}
\usepackage[T1]{fontenc}
\usepackage[utf8]{inputenc}
\usepackage{lmodern}
\usepackage{amsmath,amssymb,amsthm,mathtools}
\usepackage{mathrsfs}
\usepackage{enumitem}
\usepackage{hyperref}
\usepackage[nameinlink,capitalize]{cleveref}
\usepackage{microtype}
\usepackage{needspace}
\hypersetup{
  colorlinks=true,linkcolor=blue,citecolor=blue,urlcolor=blue,
  pdftitle={Independence of Hasse-Weil L-functions over rings of Dirichlet series},
  pdfauthor={Sho Yoshikawa}
}
\numberwithin{equation}{section}

\newtheorem{theorem}{Theorem}[section]
\newtheorem{proposition}[theorem]{Proposition}
\newtheorem{corollary}[theorem]{Corollary}
\newtheorem{lemma}[theorem]{Lemma}

\theoremstyle{definition}

\theoremstyle{remark}
\newtheorem{remark}[theorem]{Remark}

\newcommand{\QQ}{\mathbb Q}
\newcommand{\ZZ}{\mathbb Z}
\newcommand{\CC}{\mathbb C}
\newcommand{\RR}{\mathbb R}
\newcommand{\Gal}{\operatorname{Gal}}
\newcommand{\GL}{\operatorname{GL}}
\newcommand{\Tr}{\operatorname{Tr}}
\newcommand{\Frob}{\operatorname{Frob}}
\newcommand{\ud}{\operatorname{ud}}
\newcommand{\calF}{\mathcal F}
\newcommand{\calS}{\mathcal S}

\newcommand{\NT}{\mathbb N_T}

\newcommand{\Dabs}{\mathcal D_{\mathrm{abs}}}
\newcommand{\RP}{\mathcal R_P}
\newcommand{\RI}{\mathcal R_{\mathscr T}}
\newcommand{\psupp}{\operatorname{psupp}}

\title{Independence of Hasse--Weil $L$-functions\\over rings of Dirichlet series}
\author{Sho Yoshikawa\thanks{Department of Mathematics, Faculty of Science Division I, Tokyo University of Science, 1-3 Kagurazaka, Shinjuku-ku, Tokyo 162-8601, Japan.\newline
\textit{Email address:} \href{mailto:yoshikawa@rs.tus.ac.jp}{\texttt{yoshikawa@rs.tus.ac.jp}}.}}
\date{November 9, 2026}

\begin{document}
\maketitle
\begingroup
\renewcommand{\thefootnote}{}
\begin{NoHyper}
\end{NoHyper}
\endgroup

\begin{abstract}
We study linear and algebraic independence of Hasse--Weil $L$-functions of abelian varieties over $\QQ$, with Dirichlet series as coefficients. For nonzero varieties that are pairwise non-isogenous over $\overline\QQ$, we prove that their $L$-functions with arbitrary Dirichlet twists, together with all their derivatives and $1$, are linearly independent when the coefficient series use only primes from a set of natural density zero. For untwisted $L$-functions, this independence is equivalent to the condition that the varieties are pairwise non-isogenous over $\QQ$. If the individual monodromy groups are connected, the result for varieties that are pairwise non-isogenous over $\overline\QQ$ extends to larger coefficient rings: the prime supports of coefficients are unrestricted outside a prescribed set of positive upper natural density and meet it in a set of natural density zero. We also characterize algebraic independence of untwisted $L$-functions in terms of the multiplicities of the simple isogeny factors of the varieties. The proofs combine an arithmetic-function lemma of Kaczorowski--Molteni--Perelli with Frobenius density arguments and require no analytic continuation or automorphy.
\end{abstract}

\section{Introduction}

Independence questions for $L$-functions ask whether they satisfy linear or polynomial relations with coefficients in a prescribed ring of functions. Kaczorowski--Molteni--Perelli proved that distinct elements of the Selberg class $\calS$ \cite{Selberg} are linearly independent over the ring $\calF$ of $p$-finite Dirichlet series, that is, Dirichlet series that converge absolutely in a right half-plane and whose coefficients are supported on integers generated by a finite set of primes \cite{KMP}. Related results allow derivatives, meromorphic coefficients of restricted growth, and algebraic relations \cite{ChiangFeng,CimpoeasNicolae,KMPForum,Lu,Molteni,Nagoshi,Nicolae}, and Shapiro--Sparer and Pong developed an algebraic theory of arithmetic functions under Dirichlet convolution \cite{Pong,ShapiroSparer}. We compare these results with ours in \cref{rem:intro-comparison}.

In this paper we study linear and algebraic independence of Hasse--Weil $L$-functions of abelian varieties over $\QQ$, including their Dirichlet twists. By Faltings's isogeny theorem, the Hasse--Weil $L$-function determines the $\QQ$-isogeny class of the abelian variety. We ask when the assumption that the varieties are pairwise non-isogenous also rules out relations among these functions and their derivatives with Dirichlet series as coefficients. For abelian varieties over $\QQ$, we use \emph{geometrically isogenous} to mean isogenous over $\overline\QQ$.

We first prove two linear independence theorems, both including all derivatives, over different rings of Dirichlet series. We then use linear independence for products of abelian varieties to characterize algebraic independence of untwisted $L$-functions. To specify the coefficient rings, let $\Dabs$ be the ring of Dirichlet series that converge absolutely in some right half-plane, and define the prime support of $C(s)=\sum_n c(n)n^{-s}$ by
\[
 \psupp(C)=\{p:\text{$p\mid n$ for some $n$ with $c(n)\ne0$}\}.
\]
\Needspace{5\baselineskip}
Let $\pi(x)$ be the number of rational primes at most $x$. For a set $T$ of primes, its upper natural density is
\[
 \ud(T)=\limsup_{x\longrightarrow\infty}
 \frac{\#\{p\in T:p\le x\}}{\pi(x)}.
\]
If the ratio has a limit, we denote it by $\delta(T)$ and call it the natural density of $T$. We use \emph{upper density} to mean upper natural density. Our first coefficient ring is
\[
 \calF'=\{C\in\Dabs:\delta(\psupp(C))=0\}.
\]
Thus $\calF'$ contains the $p$-finite ring $\calF$: its elements may use infinitely many primes, provided those primes form a set of density zero. For a family $L_1,\ldots,L_N$ of $L$-functions, our linear independence statements assert that a relation
\begin{equation}
\label{eq:intro-relation}
 C_0(s)+\sum_{i=1}^N\sum_{k=0}^m C_{ik}(s)L_i^{(k)}(s)=0
\end{equation}
with coefficients in the specified ring forces all coefficients to vanish. Including $1$ in the family accounts for $C_0$ and lets us treat the constant term of a polynomial relation in the same way. All identities are understood in a common right half-plane of absolute convergence. The superscript $(k)$ denotes the $k$-th derivative with respect to $s$, with $L^{(0)}=L$.

Throughout the paper, Dirichlet twists are taken coefficientwise. If $L(s,A)=\sum_{n\ge1}a_A(n)n^{-s}$ is the Hasse--Weil $L$-function of $A$ and $\chi$ is a Dirichlet character, we set
\[
 L(s,A\otimes\chi)=\sum_{n\ge1}a_A(n)\chi(n)n^{-s}.
\]
Here $\chi(n)=0$ when $n$ is not coprime to the modulus of $\chi$. The twist defined by the Galois representation gives the same Euler factors outside finitely many primes; see \cref{rem:twist-convention}.

\begin{theorem}
\label{thm:main}
Let $A_1,\ldots,A_N$ be nonzero abelian varieties over $\QQ$.
\begin{enumerate}[label=\textup{(\arabic*)}]
\item Suppose that the $A_i$ are pairwise non-isogenous over $\overline\QQ$. For any Dirichlet characters $\chi_1,\ldots,\chi_N$ and every $m\ge0$, the indexed family
\[
 1,\qquad L^{(k)}(s,A_i\otimes\chi_i)
 \quad(1\le i\le N,\ 0\le k\le m)
\]
is linearly independent over $\calF'$.
\item The $A_i$ are pairwise non-isogenous over $\QQ$ if and only if
\[
 1,L(s,A_1),\ldots,L(s,A_N)
\]
are linearly independent over $\calF'$. These conditions are also equivalent to linear independence over $\calF'$ of
\[
 1,\qquad L^{(k)}(s,A_i)
 \quad(1\le i\le N,\ 0\le k\le m)
\]
for every $m\ge0$.
\end{enumerate}
\end{theorem}

Both parts of \cref{thm:main} follow from \cref{thm:representation-criterion}, which gives a necessary and sufficient condition for independence in terms of the twisted Tate modules.

Our second linear independence theorem concerns larger coefficient rings, under a connectedness assumption on the monodromy groups. For a fixed set $P$ of primes with $\ud(P)>0$, write $P^c$ for its complement and put
\[
 \RP=\{C\in\Dabs:\delta(P\cap\psupp(C))=0\}.
\]
Thus the prime support of an element $C\in\RP$ may contain any primes in $P^c$, but only a set of density zero in $P$. Equivalently, it is contained in $P^c\cup Z$ for some $Z\subseteq P$ with $\delta(Z)=0$. All densities are relative to all rational primes. The inclusion $\calF'\subseteq\RP$ is strict if and only if $\ud(P^c)>0$. For a prime $\ell$, let $T_\ell A$ be the Tate module of $A$ and write $V_\ell A=T_\ell A\otimes_{\ZZ_\ell}\QQ_\ell$. The algebraic monodromy group of $V_\ell A$ is the Zariski closure of the image of the Galois action on $V_\ell A$.

\Needspace{14\baselineskip}
\begin{theorem}
\label{thm:connected}
Let $A_1,\ldots,A_N$ be nonzero abelian varieties over $\QQ$. Suppose that, for some prime $\ell$, the algebraic monodromy group of each $V_\ell A_i$ is connected. Fix Dirichlet characters $\chi_1,\ldots,\chi_N$ and write $L_i(s)=L(s,A_i\otimes\chi_i)$. The following conditions are equivalent:
\begin{enumerate}[label=\textup{(\alph*)}]
\item The $A_i$ are pairwise non-isogenous over $\overline\QQ$.
\item For every set $P$ of primes with $\ud(P)>0$, the family $1,L_1,\ldots,L_N$ is linearly independent over $\RP$.
\item For every such $P$ and every $m\ge0$, the indexed family
\[
 1,\qquad L_i^{(k)}\quad(1\le i\le N,\ 0\le k\le m)
\]
is linearly independent over $\RP$.
\end{enumerate}
\end{theorem}

The varieties in \cref{thm:connected} need not have the same dimension. The unrestricted part of these coefficient rings can have density arbitrarily close to one. Indeed, take $P=\{p:p\equiv1\pmod q\}$ with $q\ge3$. Coefficients may then use every prime outside this progression, together with a density-zero set inside it. The unrestricted part has density $1-1/\varphi(q)$. Each independence statement concerns one fixed ring $\RP$; see \cref{rem:fixed-P}. The role of the quantifier over $P$ and the connectedness hypothesis is illustrated in \cref{rem:quadratic-twists}. Without connectedness of the monodromy groups, \cref{cor:splitting} still gives independence over $\RP$ for geometrically non-isogenous varieties, provided the intersection of $P$ with the primes splitting completely in a suitable finite Galois extension of $\QQ$ has positive upper density.

The linear independence results also lead to a criterion for algebraic independence of untwisted $L$-functions. Since
\[
 L(s,A\times B)=L(s,A)L(s,B),
\]
the assumption that the varieties are pairwise non-isogenous does not exclude polynomial relations. The precise criterion is expressed by the multiplicities of the simple isogeny factors. For a prime $\ell$, let $K_0(\QQ_\ell)$ be the Grothendieck group of finite-dimensional continuous semisimple $\QQ_\ell$-representations of $G_\QQ=\Gal(\overline\QQ/\QQ)$. In this group, $[V\oplus W]=[V]+[W]$.

\Needspace{20\baselineskip}
\begin{theorem}
\label{thm:algebraic}
Let $A_1,\ldots,A_r$ be abelian varieties over $\QQ$, fix a prime $\ell$, and write
\[
 A_i\sim_\QQ\prod_{j=1}^t B_j^{n_{ij}},
\]
where $\sim_\QQ$ denotes isogeny over $\QQ$, and the $B_j$ are $\QQ$-simple and pairwise non-isogenous over $\QQ$. The following conditions are equivalent:
\begin{enumerate}[label=\textup{(\alph*)}]
\item $L(s,A_1),\ldots,L(s,A_r)$ are algebraically independent over $\calF'$.
\item They are algebraically independent over $\CC$.
\item The vectors $(n_{ij})_{j=1}^t$ $(1\le i\le r)$ are linearly independent over $\QQ$.
\item The classes $[V_\ell A_1],\ldots,[V_\ell A_r]$ are linearly independent in $K_0(\QQ_\ell)\otimes_\ZZ\QQ$.
\end{enumerate}
If the monodromy group of $\bigoplus_i V_\ell A_i$ is connected, these conditions are also equivalent to algebraic independence over $\RP$ for every set $P$ with $\ud(P)>0$.
\end{theorem}

In particular, the $L$-functions of pairwise non-isogenous $\QQ$-simple abelian varieties are algebraically independent over $\calF'$. The assertion over $\RP$ in \cref{thm:algebraic} is proved under the assumption that the monodromy group of the direct sum $\bigoplus_i V_\ell A_i$ is connected. This assumption is stronger than the condition on the individual groups in \cref{thm:connected}. A corresponding criterion for twisted $L$-functions is given in \cref{prop:twisted-algebraic}. \Cref{rem:trace-reduction} gives an equivalent criterion in terms of Frobenius traces modulo $p$.

\Needspace{5\baselineskip}
\begin{remark}
\label{rem:intro-comparison}
We briefly compare our results with previous work.

\emph{Constant and $p$-finite coefficients.}
Kaczorowski--Molteni--Perelli proved linear independence over $\calF$ for distinct elements of $\calS$ \cite{KMP}. They also proved linear independence over $\CC$ of distinct nonconstant functions, together with all their derivatives, in a larger axiomatic class whose axioms include meromorphic continuation of finite order and a functional equation \cite{KMPForum}. For pairwise non-isogenous elliptic curves $E/\QQ$, modularity \cite{BCDT,TaylorWiles,Wiles} identifies the normalized functions $L(s+1/2,E)$ with distinct elements of $\calS$, so both independence results apply. For general abelian varieties over $\QQ$, neither automorphy nor the analytic continuation and functional equation required in these settings are known. For algebraic relations, Kaczorowski--Molteni--Perelli showed that unique factorization in $\calS$ is equivalent to algebraic independence of its distinct primitive elements over $\calF$ \cite[Corollary~2]{KMP}. Molteni proved that every nontrivial element of $\calS$ is transcendental over $\calF$ and that two distinct primitive elements are algebraically independent over $\calF$ \cite[Corollaries~2 and~3]{Molteni}. For a fixed finite Galois extension $K/\QQ$, the Artin $L$-functions attached to distinct irreducible characters of $\Gal(K/\QQ)$ are algebraically independent over $\CC$ \cite[Corollary~4]{Nicolae}.

\emph{Meromorphic coefficients.}
A lemma of Chiang--Feng \cite[Lemma~4]{ChiangFeng} implies that linear independence over $\CC$ of Dirichlet series convergent in a half-plane persists over the field of meromorphic functions $\phi$ on $\CC$ with $T(r,\phi)=o(r)$, where $T$ is the Nevanlinna characteristic. Nagoshi gives the corresponding statement for algebraic independence \cite[Proposition~4.2]{Nagoshi}, noting that it was essentially obtained by Chiang and Feng. He applies it to real shifts and derivatives of the $L$-function of an irreducible cuspidal automorphic representation of $\GL_m(\mathbb A_\QQ)$ with unitary central character \cite[Theorem~1.6]{Nagoshi}. L\"u extends these statements to meromorphic coefficients satisfying only the weaker condition $T(r_n,\phi)=o(r_n)$ along a fixed sequence $r_n\to\infty$ \cite[Theorem~1, Remark~1, and Proposition~1]{Lu}. For distinct elements of $\calS$, he proves linear independence over the ring of finite sums of products of such meromorphic functions with $p$-finite Dirichlet series \cite[Theorem~2]{Lu}. For the Artin $L$-functions above, Cimpoea\c{s}--Nicolae prove algebraic independence over the field of meromorphic functions of order less than one \cite[Corollary~9]{CimpoeasNicolae}, and L\"u allows growth along a sequence as above \cite[Corollary~2]{Lu}.

\emph{Arithmetic functions.}
Ostrowski proved that a series in $\Dabs$ with infinite prime support satisfies no algebraic differential equation over $\CC(s)$ \cite{Ostrowski}; see also \cite[Corollary~1.5]{Nagoshi}. Shapiro--Sparer developed an algebraic theory of arithmetic functions under Dirichlet convolution, using derivations attached to primes and a Jacobian criterion for algebraic independence \cite{ShapiroSparer}. Pong further developed this approach using a theorem of Ax \cite[\S3]{Pong}. These methods establish algebraic independence over rings of formal Dirichlet series annihilated by the derivations attached to a prescribed set of primes, under explicit conditions on the given arithmetic functions.

\emph{The present paper.}
We consider Hasse--Weil $L$-functions of arbitrary abelian varieties over $\QQ$ and their Dirichlet twists, using them only in a half-plane of absolute convergence. Our coefficient rings $\calF'$ and $\RP$ are defined by density conditions on prime supports. Their elements may involve infinitely many primes and need not continue meromorphically to $\CC$. Over these rings we prove linear independence of several $L$-functions together with all their derivatives, and we give criteria for linear and algebraic independence in terms of the abelian varieties and their Galois representations; see \cref{thm:representation-criterion,thm:connected,thm:algebraic}. Under the hypotheses of our independence results, Faltings's isogeny theorem, the Chebotarev density theorem, and Rajan's theorem, together with the degrees of the local polynomials, provide sets of primes on which local factors differ and which cannot be contained in the permitted coefficient supports.
\end{remark}

We now explain the ideas of the proofs. Both linear independence theorems follow from the criterion in terms of local Euler factors in \cref{prop:common}. In a proposed nontrivial relation, let $T$ be the union of the prime supports of the coefficients. The partial Euler products over $T$ are units of the coefficient ring, so the Leibniz rule gives an equivalent nontrivial relation among the remaining products and their derivatives. Comparing Dirichlet coefficients then reduces it to linear relations over $\CC$ among arithmetic functions, to which the lemma of Kaczorowski--Molteni--Perelli applies. It suffices that, outside each permitted support $T$, each pair of Euler products differs at infinitely many primes and each product has infinitely many factors different from $1$. If this condition fails, enlarging $T$ by finitely many primes yields an explicit relation from partial Euler products.

The degrees of the good local polynomials distinguish different dimensions; they also show that each Euler product has a factor different from $1$ at all but finitely many primes. For the independence assertions of \cref{thm:main}, Faltings's theorem shows that the relevant twisted Tate modules are pairwise non-isomorphic. The ordinary Chebotarev density theorem then gives, for each pair, a set of primes of positive density on which their Frobenius traces differ. A coefficient support of density zero cannot contain all but finitely many of these primes. For \cref{thm:connected}, assume that the varieties are pairwise non-isogenous over $\overline\QQ$. For a pair of equal dimension, Faltings's theorem shows that their Tate modules do not become isomorphic over any finite extension. The connectedness assumption then allows us to apply Rajan's strong multiplicity one theorem and the algebraic Chebotarev density theorem, as in \cref{lem:Rajan-density}: the set of primes with equal twisted traces has natural density zero. Thus, for every $P$ with $\ud(P)>0$ and every permitted support $T$, infinitely many primes of $P\setminus T$ distinguish each pair of Euler products. Conversely, isomorphic twisted Tate modules give a relation over $\calF'$. If two varieties are geometrically isogenous, their twisted local factors agree at almost all primes splitting completely in a finite Galois extension over which an isogeny is defined and the characters are trivial. Taking $P$ to be this set of split primes gives a relation over $\RP$.

For algebraic independence, we prove linear independence of the $L$-functions of the products $A_\alpha=\prod_i A_i^{\alpha_i}$ for $\alpha\in\ZZ_{\ge0}^r\setminus\{0\}$, together with $1$. The identity
\[
 \prod_i L(s,A_i)^{\alpha_i}
 =L(s,A_\alpha)
\]
turns a polynomial relation among the $L(s,A_i)$ into a linear relation among these $L(s,A_\alpha)$. Under condition \textup{(d)} of \cref{thm:algebraic}, the Tate modules of the $A_\alpha$ are nonzero and pairwise non-isomorphic, so \cref{thm:representation-criterion}, applied to any finite collection of these products, excludes such relations. Conversely, dependent multiplicity vectors give an equality between two distinct monomials. For the assertion over $\RP$, connectedness of the monodromy group of $\bigoplus_i V_\ell A_i$ passes to the direct sum of the Tate modules of any finite collection of the $A_\alpha$. By \cref{lem:joint-density}, the algebraic Chebotarev density theorem then gives natural density zero for the primes with equal traces, and \cref{cor:density-criteria} yields the required linear independence.

The paper is organized as follows. \Cref{sec:linear} develops the coefficient rings and proves the common independence criterion. \Cref{sec:Galois} supplies the Frobenius density estimates and fixes the conventions for local factors and twists. \Cref{sec:applications} proves the linear independence theorems and derives the algebraic independence criteria; \cref{cor:partial-products,cor:curves} give applications to partial Euler products and curves. \Cref{rem:curve-hypotheses} provides sufficient conditions for the hypotheses concerning Jacobians.

\paragraph{Acknowledgements.}
The author thanks Yuta Katayama for drawing his attention to papers on the independence of $L$-functions, which led to the present work.

\paragraph{Disclosure on the use of generative AI.}
The ideas and proofs presented in this paper are the author’s own, and the initial manuscript was written by the author. Generative AI tools, specifically ChatGPT (Sol 6.1) and Claude (5.5), were used to assist in organizing and refining some of the proofs, improving the exposition and language, and conducting bibliographic searches. The author independently verified all mathematical arguments, statements, and references. The author takes full responsibility for the content of this paper.

\section{Coefficient rings and a common independence criterion}
\label{sec:linear}

\subsection{Prime supports and coefficient rings}

For a set $T$ of primes, write
\[
 \NT=\{n\ge1:\text{every prime divisor of $n$ belongs to $T$}\},
\]
\[
 \mathcal D_T=\left\{\sum_{n\ge1}c(n)n^{-s}\in\Dabs:
 c(n)=0\text{ for }n\notin\NT\right\}.
\]
Note that $1\in\NT$ even if $T=\varnothing$. A prime support for a Dirichlet series $C$ is a set $T$ of primes such that $C\in\mathcal D_T$. Its smallest prime support is $\psupp(C)$.

A family $\mathscr T$ of prime sets is called \emph{admissible} if it is closed under taking subsets and finite unions and contains every finite set. Fix an admissible family $\mathscr T$ throughout this section, except when another family is specified, and define
\[
 \RI=\bigcup_{T\in\mathscr T}\mathcal D_T.
\]
For the two coefficient rings introduced in the introduction, set
\[
 \mathscr T_0=\{T:\delta(T)=0\},\qquad
 \mathscr T_P=\{T:\delta(T\cap P)=0\},
\]
where $P$ is a fixed set of primes with $\ud(P)>0$. Both families are admissible, and
\[
 \calF'=\mathcal R_{\mathscr T_0},\qquad
 \RP=\mathcal R_{\mathscr T_P}.
\]
Thus $\calF'$ allows only a density-zero set of supporting primes, while $\RP$ allows arbitrary supports in $P^c$ and density-zero supports in $P$. The density of $T\cap P$ is taken relative to all rational primes. We identify two Dirichlet series when they agree in a common right half-plane. By the uniqueness theorem for Dirichlet series \cite[Theorem~11.3]{Apostol}, this is equivalent to equality of their coefficients.

\begin{lemma}
\label{lem:ring}
The ring $\RI$ is a differential $\CC$-algebra and an integral domain. If $C(s)=\sum_n c(n)n^{-s}\in\mathcal D_T$, $T\in\mathscr T$, and $c(1)\ne0$, then $C$ is a unit in $\RI$ and $C^{-1}\in\mathcal D_T$.
\end{lemma}

\begin{proof}
Constants belong to $\RI$. Addition and Dirichlet convolution preserve support in the union of the prime supports, which still belongs to $\mathscr T$. Absolute convergence is preserved in a common right half-plane. Termwise differentiation gives
\[
 C^{(j)}(s)=\sum_{n\geq1}c(n)(-\log n)^j n^{-s}.
\]
The resulting series is still supported on $\NT$. It converges absolutely sufficiently far to the right: for every $\varepsilon>0$ there is a constant $M_{j,\varepsilon}$, depending only on $j$ and $\varepsilon$, such that $(\log n)^j\le M_{j,\varepsilon}n^\varepsilon$ for $n\ge1$.

Let $C(s)=\sum_n c(n)n^{-s}$ and $D(s)=\sum_n d(n)n^{-s}$ be nonzero Dirichlet series, and let $n_0$ and $m_0$ be the smallest indices with $c(n_0)\ne0$ and $d(m_0)\ne0$, respectively. In the coefficient of $(n_0m_0)^{-s}$ in $CD$, a nonzero term $c(n)d(m)$ must satisfy $n\ge n_0$, $m\ge m_0$, and $nm=n_0m_0$. Hence the coefficient is $c(n_0)d(m_0)\ne0$, so there are no zero divisors.

To construct the inverse, write
\[
 C(s)=c(1)(1+H(s)),\qquad
 H(s)=\sum_{n\ge2}\frac{c(n)}{c(1)}n^{-s}.
\]
Thus $H$ has no term with index $n=1$. For a Dirichlet series $D(s)=\sum_n d(n)n^{-s}$ and $\sigma\in\RR$ at which $D$ converges absolutely, put
\[
 \|D\|_\sigma=\sum_{n\ge1}|d(n)|n^{-\sigma}.
\]
This is the $\ell^1$ norm of the sequence $(d(n)n^{-\sigma})_{n\ge1}$. If $D_1$ and $D_2$ both converge absolutely at $\sigma$, Dirichlet convolution gives $\|D_1D_2\|_\sigma\le\|D_1\|_\sigma\|D_2\|_\sigma$. Since the series defining $H$ starts at $n=2$, dominated convergence gives $\|H\|_\sigma\to0$ as $\sigma\to+\infty$. Choose $\sigma$ with $\|H\|_\sigma<1$. Then $\sum_{q\ge0}\|H^q\|_\sigma<\infty$, so
\[
 c(1)^{-1}\sum_{q\ge0}(-H)^q
\]
defines a Dirichlet series supported on $\NT$ that converges absolutely at $\sigma$. Multiplying it by $C=c(1)(1+H)$ gives $1$, as required.
\end{proof}

Taking $\mathscr T$ to be the family of all sets of primes gives $\RI=\Dabs$, so the lemma also shows that $\Dabs$ is an integral domain.

An element of the fraction field of a coefficient ring is a quotient of two Dirichlet series and is meromorphic on some right half-plane. Multiplying a finite relation over the fraction field by a common nonzero denominator gives a relation over the ring itself. Thus linear and algebraic independence over a coefficient ring are equivalent to the corresponding independence over its fraction field.

If $\ud(P^c)>0$, then $\sum_{p\notin P}p^{-s}$ belongs to $\RP$ but not to $\calF'$. Conversely, if $\ud(P^c)=0$, then $\delta(T\cap P)=0$ implies $\delta(T)=0$, so $\RP=\calF'$. Neither coefficient ring imposes global meromorphic continuation or a growth condition on its elements.

\subsection{Separating Euler factors}

Let $\calS'$ be the class of series $L(s)=\sum_n a(n)n^{-s}\in\Dabs$ with multiplicative coefficients and $a(1)=1$. Here multiplicative means $a(nr)=a(n)a(r)$ for coprime $n,r$. For any set $T$ of primes, define
\[
 E_T(s)=\sum_{n\in\NT}a(n)n^{-s},\qquad
 L^{(T)}(s)=\sum_{\substack{n\ge1\\p\mid n\Rightarrow p\notin T}}a(n)n^{-s}.
\]
The series $E_T$ is the product of the Euler factors at primes in $T$, and $L^{(T)}$ is the product at the remaining primes. A superscript $(k)$ with $k$ a nonnegative integer denotes a derivative.

\Needspace{14\baselineskip}
\begin{lemma}
\label{lem:depletion}
For $L\in\calS'$ and any $T$,
\begin{equation}
\label{eq:T-factorization}
 L=E_TL^{(T)}
\end{equation}
in a right half-plane of absolute convergence. If $T\in\mathscr T$, then $E_T$ is a unit in $\RI$, and for every $m\ge0$ the identities
\begin{equation}
\label{eq:Leibniz}
 L^{(k)}=\sum_{h=0}^k\binom{k}{h}E_T^{(k-h)}
 \bigl(L^{(T)}\bigr)^{(h)}\quad(0\le k\le m)
\end{equation}
express
\[
 (L,L',\ldots,L^{(m)})
 \quad\text{in terms of}\quad
 \bigl(L^{(T)},(L^{(T)})',\ldots,(L^{(T)})^{(m)}\bigr)
\]
by an invertible triangular matrix over $\RI$, with $E_T$ on every diagonal entry.
\end{lemma}

\begin{proof}
Every integer has a unique factorization $n=dr$ with $d\in\NT$ and no prime divisor of $r$ in $T$. Then $(d,r)=1$ and $a(n)=a(d)a(r)$, proving \eqref{eq:T-factorization}. Both factors are subseries of an absolutely convergent series. If $T\in\mathscr T$, then $E_T\in\mathcal D_T$ has coefficient $1$ at $1$, so \cref{lem:ring} applies. Formula \eqref{eq:Leibniz} is the Leibniz rule, and its triangular matrix has the unit $E_T$ on the diagonal.
\end{proof}

\subsection{An independence criterion}

For a multiplicative arithmetic function $f$, put
\[
 f^{[k]}(n)=(-\log n)^kf(n)\quad(k\ge0),
\]
with $(-\log n)^0=1$. Two multiplicative functions $f,g$ are equivalent if, outside a finite set of primes, $f(p^j)=g(p^j)$ for every $j\ge1$. Let $e$ be the arithmetic function with $e(1)=1$ and $e(n)=0$ for $n>1$. It is the coefficient function of the constant Dirichlet series $1$.

\begin{lemma}[Kaczorowski--Molteni--Perelli]
\label{lem:KMP}
If $e,f_1,\ldots,f_N$ are pairwise non-equivalent multiplicative arithmetic functions, then for every $m\ge0$ the indexed family
\[
 e,\qquad f_i^{[k]}\quad(1\le i\le N,\ 0\le k\le m)
\]
is linearly independent over $\CC$.
\end{lemma}

\begin{proof}
The assertion without $e$ is \cite[Lemma~1]{KMPForum}; see also \cite[Lemma~3]{CimpoeasNicolae}. Suppose that
\[
 c_0 e(n)+\sum_{i=1}^N\sum_{k=0}^m c_{ik}f_i^{[k]}(n)=0
 \qquad(n\ge1)
\]
for complex constants $c_0,c_{ik}$. For each $n$, multiply this equality by $-\log n$. The first term becomes zero: $\log 1=0$ and $e(n)=0$ for $n>1$. The other terms give
\[
 \sum_{i=1}^N\sum_{k=0}^m c_{ik}f_i^{[k+1]}(n)=0
 \qquad(n\ge1).
\]
The cited lemma, applied with derivatives up to order $m+1$, gives $c_{ik}=0$ for every $i,k$. The original equality at $n=1$ then gives $c_0=0$.
\end{proof}

For $L_i(s)=\sum_n a_i(n)n^{-s}\in\calS'$, define the local factors as formal power series
\[
 L_{i,p}(X)=\sum_{\nu\ge0}a_i(p^\nu)X^\nu,
\]
and put
\begin{equation}
\label{eq:local-difference-sets}
 D_i=\{p:L_{i,p}\ne1\},\qquad
 D_{ij}=\{p:L_{i,p}\ne L_{j,p}\}\quad(i\ne j).
\end{equation}
Two local factors at $p$ are equal exactly when their coefficients at every power of $X$ agree. For each $i$, write $E_{i,T}$ and $L_i^{(T)}$ for the two series obtained from $L_i$ as in \cref{lem:depletion}.

\Needspace{13\baselineskip}
\begin{proposition}
\label{prop:common}
Let $\mathscr T$ be an admissible family of prime sets and let $L_1,\ldots,L_N\in\calS'$. The following conditions are equivalent:
\begin{enumerate}[label=\textup{(\alph*)}]
\item For every $T\in\mathscr T$, the set $D_i\setminus T$ is infinite for every $i$, and $D_{ij}\setminus T$ is infinite for every pair $i\ne j$.
\item The family $1,L_1,\ldots,L_N$ is linearly independent over $\RI$.
\item For every $m\ge0$, the indexed family
\[
 1,\qquad L_i^{(k)}\quad(1\le i\le N,\ 0\le k\le m)
\]
is linearly independent over $\RI$.
\end{enumerate}
\end{proposition}

\begin{proof}
We first prove \textup{(a)}$\Rightarrow$\textup{(c)}. For every $T\in\mathscr T$, condition \textup{(a)} says exactly that
\begin{equation}
\label{eq:non-equivalence}
 e,a_1^{(T)},\ldots,a_N^{(T)}
 \quad\text{are pairwise non-equivalent},
\end{equation}
where $a_i^{(T)}$ is the coefficient function of $L_i^{(T)}$.
Suppose that a nontrivial relation
\begin{equation}
\label{eq:relation}
 C_0+\sum_{i=1}^N\sum_{k=0}^m C_{ik}L_i^{(k)}=0,
 \qquad C_0,C_{ik}\in\RI,
\end{equation}
holds. Its coefficients have prime supports contained in a common $T\in\mathscr T$. Write $L_i=E_{i,T}L_i^{(T)}$. By \cref{lem:depletion} and \eqref{eq:Leibniz}, the relation becomes
\begin{equation}
\label{eq:depleted-relation}
 C_0+\sum_{i=1}^N\sum_{h=0}^m
 B_{ih}\bigl(L_i^{(T)}\bigr)^{(h)}=0,
\end{equation}
where
\[
 B_{ih}=\sum_{k=h}^m\binom{k}{h}C_{ik}E_{i,T}^{(k-h)}.
\]
Every $B_{ih}$ lies in $\mathcal D_T$. The invertibility in \cref{lem:depletion} shows that $C_0$ and the $B_{ih}$ are not all zero.

Write $C_0=\sum_d c_0(d)d^{-s}$ and $B_{ih}=\sum_d b_{ih}(d)d^{-s}$. Every integer $n$ has a unique decomposition
\[
 n=dr,\qquad d\in\NT,\qquad p\mid r\ \Longrightarrow\ p\notin T,
\]
obtained by separating its prime factors in $T$ from those outside $T$. Fix $d\in\NT$ and let $r$ have no prime divisor in $T$. The coefficient of $(dr)^{-s}$ in $B_{ih}(L_i^{(T)})^{(h)}$ is
\[
 \sum_{uv=dr}b_{ih}(u)a_i^{(T)}(v)(-\log v)^h
 =b_{ih}(d)a_i^{(T)}(r)(-\log r)^h.
\]
Indeed, $b_{ih}(u)=0$ unless $u\in\NT$, while $a_i^{(T)}(v)=0$ if a prime in $T$ divides $v$. Thus the only term that can be nonzero has $u=d$ and $v=r$. Similarly, the coefficient of $(dr)^{-s}$ in $C_0$ is $c_0(d)e(r)$: it is zero unless $r=1$. The uniqueness theorem for Dirichlet series, applied to \eqref{eq:depleted-relation}, therefore gives
\begin{equation}
\label{eq:arithmetic-relation}
 c_0(d)e(r)+\sum_{i=1}^N\sum_{h=0}^m
 b_{ih}(d)a_i^{(T)}(r)(-\log r)^h=0.
\end{equation}
If $r$ has a prime divisor in $T$, all the arithmetic functions in this equality vanish, so it holds for every $r\ge1$. By \eqref{eq:non-equivalence} and \cref{lem:KMP}, all its coefficients are zero. Varying $d$ gives $C_0=B_{ih}=0$, a contradiction.

The implication \textup{(c)}$\Rightarrow$\textup{(b)} follows by taking $m=0$.

To prove \textup{(b)}$\Rightarrow$\textup{(a)}, suppose that \textup{(a)} fails. Choose $T\in\mathscr T$ for which some $D_i\setminus T$ or $D_{ij}\setminus T$ is finite. After adjoining these finitely many primes, we still have $T\in\mathscr T$. If $D_{ij}\subseteq T$, then $L_i^{(T)}=L_j^{(T)}$, so
\[
 E_{j,T}L_i-E_{i,T}L_j=0
\]
is a nontrivial relation over $\RI$. If $D_i\subseteq T$, then $L_i=E_{i,T}\in\mathcal D_T$, giving a relation between $1$ and $L_i$. Both alternatives contradict \textup{(b)}.
\end{proof}

\Needspace{15\baselineskip}
\begin{corollary}
\label{cor:relation-support}
Let $L_1,\ldots,L_N\in\calS'$ and fix a set $T$ of primes. Put
\[
 \mathscr T(T)=\{U:U\setminus T\text{ is finite}\},\qquad
 \mathcal R_{\mathscr T(T)}
 =\{C\in\Dabs:\psupp(C)\setminus T\text{ is finite}\}.
\]
There is a nontrivial relation over $\mathcal R_{\mathscr T(T)}$ among $1$, the $L_i$, and finitely many of their derivatives if and only if at least one of the following holds:
\begin{enumerate}[label=\textup{(\roman*)}]
\item For some $i$, the local factor $L_{i,p}$ is $1$ at all but finitely many primes $p\notin T$.
\item For some pair $i\ne j$, the local factors $L_{i,p}$ and $L_{j,p}$ agree at all but finitely many primes $p\notin T$.
\end{enumerate}
In either case there is already a nontrivial relation among $1,L_1,\ldots,L_N$ over $\mathcal R_{\mathscr T(T)}$.

In particular, given a nontrivial relation with coefficients in $\Dabs$, let $T$ be the union of the prime supports of those coefficients. Then $T$ contains all but finitely many primes in at least one of the sets $D_i$ or $D_{ij}$.
\end{corollary}

\begin{proof}
The family $\mathscr T(T)$ is admissible. Apply \cref{prop:common} with this family. Dependence is equivalent to the finiteness of some $D_i\setminus T$ or $D_{ij}\setminus T$. These are precisely conditions \textup{(i)} and \textup{(ii)}, respectively. The proposition also gives a relation involving no derivatives.
\end{proof}

\Needspace{14\baselineskip}
\begin{corollary}
\label{cor:density-criteria}
Let $L_1,\ldots,L_N\in\calS'$ and let $P$ be a set of primes with $\ud(P)>0$. The following conditions are equivalent:
\begin{enumerate}[label=\textup{(\alph*)}]
\item For every $i$, $\ud(P\cap D_i)>0$, and for every pair $i\ne j$, $\ud(P\cap D_{ij})>0$.
\item The family $1,L_1,\ldots,L_N$ is linearly independent over $\RP$.
\item For every $m\ge0$, the indexed family
\[
 1,\qquad L_i^{(k)}\quad(1\le i\le N,\ 0\le k\le m)
\]
is linearly independent over $\RP$.
\end{enumerate}
Taking $P$ to be the set of all primes gives the corresponding equivalence for $\calF'$. Also, if every $D_i$ and $D_{ij}$ has natural density one, independence holds over every $\RP$ with $\ud(P)>0$.
\end{corollary}

\begin{proof}
Apply \cref{prop:common} with $\mathscr T=\mathscr T_P$. Its conditions \textup{(b)} and \textup{(c)} are the corresponding conditions here. To identify condition \textup{(a)}, let $D$ be any one of the sets $D_i$ or $D_{ij}$. If $\ud(P\cap D)>0$, then $P\cap D$ remains infinite after deleting any $T$ with $\delta(P\cap T)=0$. Conversely, if $\ud(P\cap D)=0$, then $D\in\mathscr T_P$, so the condition that $D\setminus T$ be infinite for every $T\in\mathscr T_P$ fails. This proves the equivalence. The final assertion follows because deleting a density-zero set preserves upper density.
\end{proof}

\begin{proposition}
\label{prop:Euler-modifications}
Let $L_1,\ldots,L_N\in\Dabs$, and let $U_1,\ldots,U_N$ be units of $\RI$. For every $m\ge0$, the two indexed families
\[
 \begin{gathered}
 1,\quad L_i^{(k)}\quad(1\le i\le N,\ 0\le k\le m),\\[2pt]
 1,\quad (U_iL_i)^{(k)}\quad(1\le i\le N,\ 0\le k\le m)
 \end{gathered}
\]
are either both linearly independent over $\RI$ or both linearly dependent over $\RI$. Algebraic independence of the $L_i$ is also preserved by these multiplications. In particular, the conclusions over $\calF'$ and $\RP$ are preserved by finitely many changes of Euler factors whose local factors and their inverses converge absolutely sufficiently far to the right.
\end{proposition}

\begin{proof}
The Leibniz rule gives an invertible triangular transformation with diagonal entries $U_i$; the constant function $1$ is unchanged. For algebraic independence, the substitutions $X_i\mapsto U_iX_i$ give an automorphism of $\RI[X_1,\ldots,X_N]$. For finitely many local-factor changes, the multiplier and its inverse are Dirichlet series supported on finitely many primes, hence units of both rings.
\end{proof}

\begin{remark}
\label{rem:fixed-P}
Each independence statement over $\RP$ concerns a fixed set $P$. To see why this matters, partition the primes into sets $T_1,T_2$ of density $1/2$. A ring allowing all Dirichlet series supported on either set contains both partial Euler products $E_{T_1}$ and $E_{T_2}$ of any Hasse--Weil $L$-function $L$. It therefore contains
\[
 L=E_{T_1}E_{T_2},
\]
so $1$ and $L$ are linearly dependent over that ring. Thus one cannot allow all prime supports of density less than one simultaneously.
\end{remark}

\section{Frobenius traces and local factors}
\label{sec:Galois}

Throughout this section, $\ell$ is a prime and $F$ is an $\ell$-adic field, meaning a finite extension of $\QQ_\ell$. For a number field $K\subseteq\overline\QQ$, write $G_K=\Gal(\overline\QQ/K)$. At a prime $p$ where a representation of $G_\QQ$ is unramified, $\Frob_p$ denotes arithmetic Frobenius.

\subsection{Frobenius traces}

The algebraic monodromy group of a continuous representation $\rho:G_\QQ\to\GL(V)$ on a finite-dimensional $F$-vector space is the Zariski closure of its image.

\begin{lemma}
\label{lem:trace-density}
Let
\[
 \rho_i:G_\QQ\longrightarrow\GL_{r_i}(F)\quad(i=1,2),
\]
be continuous semisimple representations, unramified outside finitely many primes. If $\rho_1\not\simeq\rho_2$, the set of common unramified primes satisfying
\[
 \Tr\rho_1(\Frob_p)\ne\Tr\rho_2(\Frob_p)
\]
contains a subset of positive natural density.
\end{lemma}

\begin{proof}
Over a field of characteristic zero, semisimple representations with the same trace are isomorphic, by the Brauer--Nesbitt theorem. Hence
\[
 t(\sigma)=\Tr\rho_1(\sigma)-\Tr\rho_2(\sigma)
\]
is a nonzero continuous class function. Choose $\sigma_0$ with $t(\sigma_0)\ne0$. By continuity, $t$ is nonzero in an open neighborhood of $\sigma_0$. Choose an open normal subgroup $U$ of $G_\QQ$ such that $\sigma_0U$ lies in this neighborhood. Since $t$ is invariant under conjugation, it is nonzero on every conjugate of $\sigma_0U$ as well.

Let $M/\QQ$ be the finite Galois extension fixed by $U$, and let $C$ be the conjugacy class of the image of $\sigma_0$ in $\Gal(M/\QQ)$. By the Chebotarev density theorem, the primes unramified in $M$ with Frobenius class $C$ have natural density $|C|/[M:\QQ]>0$; see \cite[\S2.1, Th\'eor\`eme~1]{SerreChebotarev}. Exclude the finitely many primes where either representation is ramified. At every remaining prime in this set, $t(\Frob_p)\ne0$, so the two traces differ.
\end{proof}

For a continuous representation $\rho:G_\QQ\to\GL(V)$ on a finite-dimensional $F$-vector space with monodromy group $G$, its \emph{connectedness field} is the finite Galois extension $K_\rho/\QQ$ defined by
\[
 G_{K_\rho}=\rho^{-1}(G^0(F)).
\]
Here $G^0$ is the identity component of $G$. The image of $G_{K_\rho}$ is Zariski dense in $G^0$.

We also use the following algebraic form of the Chebotarev density theorem.

\begin{theorem}[Algebraic Chebotarev density theorem]
\label{thm:algebraic-chebotarev}
Let $\rho:G_\QQ\to\GL(V)$ be continuous and unramified outside finitely many primes, and let $G$ be its algebraic monodromy group. Let $X\subseteq G$ be a closed subvariety defined over $F$ and stable under conjugation by $G$. Let $c$ be the number of connected components of $G$ over an algebraic closure of $F$, and $c_X$ the number of these components contained in $X$. Then
\[
 \delta\bigl(\{p:\rho\text{ is unramified at }p,\ \rho(\Frob_p)\in X(F)\}\bigr)
 =\frac{c_X}{c}.
\]
\end{theorem}

This is \cite[Theorem~3]{Rajan}, also stated in \cite[Theorem~2.2]{PatankarRajan}, specialized to the base field $\QQ$. The asymptotic formula in Rajan's theorem gives natural density. His algebraic formulation builds on Serre's $\ell$-adic Chebotarev method \cite[\S4]{SerreChebotarev}.

\begin{lemma}
\label{lem:Rajan-density}
Let $\rho_1,\rho_2:G_\QQ\to\GL_r(F)$ be continuous semisimple representations, unramified outside finitely many primes. Suppose that the monodromy group of $\rho_1$ is connected and that $\rho_1,\rho_2$ do not become isomorphic over any finite extension of $\QQ$. Let $\eta_1,\eta_2:G_\QQ\to F^\times$ be continuous characters of finite order. Then the set of primes unramified for both representations and both characters at which
\[
 \Tr(\rho_1\otimes\eta_1)(\Frob_p)=\Tr(\rho_2\otimes\eta_2)(\Frob_p)
\]
has natural density zero.
\end{lemma}

\begin{proof}
Put $\tau=\rho_2\otimes(\eta_2\eta_1^{-1})$ and let $G$ be the monodromy group of $\rho_1\oplus\tau$, viewed as a subgroup of $\GL_r\times\GL_r$. The closed subvariety
\[
 X=\{(g_1,g_2)\in G:\Tr g_1=\Tr g_2\}
\]
is stable under conjugation by $G$. Outside the finitely many ramified primes, the set in the statement is
\[
 \{p:(\rho_1\oplus\tau)(\Frob_p)\in X(F)\},
\]
since dividing the original trace equality by $\eta_1(\Frob_p)$ gives $\Tr\rho_1(\Frob_p)=\Tr\tau(\Frob_p)$. By \cref{thm:algebraic-chebotarev}, this set has a natural density.

Suppose that its density is positive. Here we use connectedness of the monodromy group of $\rho_1$: it allows us to apply Rajan's strong multiplicity one theorem for $\ell$-adic representations \cite[Theorem~2(i)]{Rajan}, in the form of \cite[Theorem~2.3]{PatankarRajan}. This gives
\[
 \rho_1|_{G_M}\simeq\tau|_{G_M}
\]
for some finite extension $M/\QQ$. Over a further finite extension on which $\eta_2\eta_1^{-1}$ is trivial, the restrictions of $\rho_1$ and $\rho_2$ are therefore isomorphic. This contradicts the hypothesis, so the density is zero.
\end{proof}

\begin{lemma}
\label{lem:joint-density}
Let $\rho_i:G_\QQ\to\GL_{r_i}(F)$ $(1\le i\le N)$ be pairwise non-isomorphic continuous semisimple representations, and let $S$ be a finite set of primes containing every prime where at least one $\rho_i$ is ramified. If the monodromy group of the direct sum $\bigoplus_i\rho_i$ is connected, then
\[
 \delta\bigl(\{p\notin S:\Tr\rho_i(\Frob_p)=\Tr\rho_j(\Frob_p)
             \text{ for some }i<j\}\bigr)=0.
\]
\end{lemma}

\begin{proof}
Let $G$ be the monodromy group of $\bigoplus_i\rho_i$, viewed as a subgroup of $\prod_i\GL_{r_i}$. For each $i<j$, the function $(g_k)_k\mapsto\Tr g_i-\Tr g_j$ is a regular function on $G$ invariant under conjugation. It is not identically zero, since otherwise the two representations would have the same trace and hence would be isomorphic by Brauer--Nesbitt. Since $G$ is connected, the zero locus contains no connected component of $G$. By \cref{thm:algebraic-chebotarev}, the primes with equal traces for this pair have natural density zero. The union over all pairs has natural density zero as well.
\end{proof}

We call the monodromy group of a direct sum the \emph{joint monodromy group} of its summands. Its connectedness implies that of each summand's group, since the latter is the image under the corresponding projection. In \cref{thm:algebraic,prop:twisted-algebraic}, we use \cref{lem:joint-density} for the representations associated with products of $L$-functions. For example, for nonnegative integers $a,b$, the product $L(s,A_1)^aL(s,A_2)^b$ corresponds to the direct sum of $a$ copies of $V_\ell A_1$ and $b$ copies of $V_\ell A_2$.

\subsection{Tate modules and local factors}

Let $A$ be an abelian variety over a number field $K$. Its Tate module and rational Tate module are
\[
 T_\ell A=\varprojlim_n A[\ell^n](\overline\QQ),\qquad
 V_\ell A=T_\ell A\otimes_{\ZZ_\ell}\QQ_\ell.
\]
Write $\rho_{A,\ell}$ for the action of $G_K$ on $V_\ell A$. Faltings proved that $V_\ell A$ is semisimple as a $G_K$-representation and that for any abelian varieties $A,B$ over $K$ the natural map
\[
 \operatorname{Hom}_K(A,B)\otimes_\ZZ\QQ_\ell
 \longrightarrow\operatorname{Hom}_{G_K}(V_\ell A,V_\ell B)
\]
is an isomorphism \cite{Faltings}. In particular, $A$ and $B$ are $K$-isogenous if and only if their rational Tate modules are isomorphic to each other as $G_K$-representations.

Now let $A$ be defined over $\QQ$. For a prime $p\ne\ell$, let $I_p$ be the inertia subgroup of a decomposition group at $p$. Arithmetic Frobenius acts on the coinvariants $(V_\ell A)_{I_p}$ through the quotient of the decomposition group by $I_p$. Define the local polynomial and the Euler factor by
\begin{equation}
\label{eq:local-factor}
 P_p(A,X)=\det\!\left(1-\Frob_p X\mid(V_\ell A)_{I_p}\right),
 \qquad L_p(A,X)=P_p(A,X)^{-1}.
\end{equation}
This agrees with the definition using $H^1$ in \'etale cohomology. Indeed, there are natural isomorphisms
\[
 H^1_{\mathrm{\acute et}}(A_{\overline\QQ},\QQ_\ell)
 \simeq(V_\ell A)^\vee,\qquad
 \bigl((V_\ell A)^\vee\bigr)^{I_p}
 \simeq\bigl((V_\ell A)_{I_p}\bigr)^\vee,
\]
where $\vee$ denotes the $\QQ_\ell$-linear dual. On the dual representation, geometric Frobenius acts by the transpose of arithmetic Frobenius on the coinvariants. Hence the determinant defining the local polynomial from cohomology is $P_p(A,X)$.

The inertia invariants $(V_\ell A)^{I_p}$ identify with the rational Tate module of the identity component of the special fiber of the N\'eron model \cite[Expos\'e~IX, \S2.2, especially Proposition~2.2.5]{SGA7I}. This identity component, modulo its unipotent part, is an extension of an abelian variety $B_p$ by a torus $T_p$; see \cite[Expos\'e~IX, \S2]{SGA7I}. Since $\QQ_\ell(1)$ is unramified at $p$, a polarization identifies
\[
 (V_\ell A)_{I_p}\simeq
 \bigl((V_\ell A)^{I_p}\bigr)^\vee(1).
\]
Consequently $P_p(A,X)$ is the product of the polynomials $\det(1-\Frob_p X)$ on $X^*(T_p)\otimes_\ZZ\QQ_\ell$ and on $V_\ell(B_p^\vee)$, where $X^*(T_p)$ is the character lattice of the torus and $B_p^\vee$ is the dual abelian variety. Frobenius acts on the integral lattice $X^*(T_p)$. The Frobenius endomorphism of $B_p^\vee$ has an integral characteristic polynomial independent of $\ell$ by \cite[\S19, Theorem~4]{Mumford}. Thus $P_p(A,X)\in1+X\ZZ[X]$ is independent of the choice of $\ell\ne p$. For each prime $p$, we use any auxiliary prime $\ell\ne p$ in its definition.

We put
\[
 L(s,A)=\prod_p L_p(A,p^{-s})=\sum_{n\ge1}a_A(n)n^{-s}.
\]
The coefficients are multiplicative and $a_A(1)=1$. At a prime of good reduction with $p\ne\ell$,
\begin{equation}
\label{eq:prime-trace}
 P_p(A,X)=\det(1-\rho_{A,\ell}(\Frob_p)X\mid V_\ell A),
 \qquad a_A(p)=\Tr\rho_{A,\ell}(\Frob_p).
\end{equation}
Writing $g=\dim A$, the coefficient of $X^{2g}$ in this polynomial is
\[
 \det\!\left(-\rho_{A,\ell}(\Frob_p)\mid V_\ell A\right)\ne0,
\]
since Frobenius acts invertibly. Thus $\deg P_p(A,X)=2g$. For an elliptic curve $A$, the polynomial is $1-a_p(A)X+pX^2$, where $a_p(A)=a_A(p)$. The Weil bounds imply absolute convergence of the product over good primes for $\Re(s)>3/2$. The finitely many remaining factors converge sufficiently far to the right, so $L(s,A)\in\calS'$.

\subsection{Dirichlet twists}

Let $\chi$ be a Dirichlet character modulo $q$, extended by zero on integers not coprime to $q$. Define
\begin{equation}
\label{eq:twist-definition}
 L(s,A\otimes\chi)=\sum_{n\ge1}a_A(n)\chi(n)n^{-s}
                 =\prod_p L_p(A,\chi(p)p^{-s}).
\end{equation}
The coefficients are multiplicative and the series converges absolutely wherever the untwisted series does. Let $E\subset\CC$ be a number field containing the values of $\chi$ and choose an embedding $\iota_\ell:E\hookrightarrow F$, enlarging the $\ell$-adic field $F$ if necessary. The finite-order character $\chi_\ell:G_\QQ\to F^\times$ associated with the primitive character inducing $\chi$ satisfies
\[
 \chi_\ell(\Frob_p)=\iota_\ell(\chi(p))
 \qquad(p\nmid q\ell).
\]
Here and below we extend scalars on the Tate module to $F$ when forming its twist by $\chi_\ell$.

If $A$ has good reduction at $p$ and $p\nmid q\ell$, applying $\iota_\ell$ to the coefficients of $P_p(A,\chi(p)X)$ gives
\[
 \det\!\left(1-(\rho_{A,\ell}\otimes\chi_\ell)(\Frob_p)X\right).
\]
In particular,
\begin{equation}
\label{eq:twisted-prime-trace}
 \iota_\ell(a_A(p)\chi(p))
 =\Tr\!\left((\rho_{A,\ell}\otimes\chi_\ell)(\Frob_p)\right).
\end{equation}
The condition $p\nmid q$ also excludes primes dividing the modulus but not the conductor of an imprimitive character. The polynomial has degree $2\dim A$, so for $A\ne0$ the twisted Euler factor differs from $1$ at all but finitely many primes. Since $\iota_\ell$ is injective, two local polynomials over $E$ agree if and only if their images in $F[X]$ agree. Scalar extension in characteristic zero and twisting by a character preserve semisimplicity, so the twisted Tate module is semisimple.

\begin{remark}
\label{rem:twist-convention}
The twisted Tate module $W=(V_\ell A\otimes_{\QQ_\ell}F)\otimes_F\chi_\ell$ gives the local factor
\[
 \det\!\left(1-\Frob_p X\mid W_{I_p}\right)^{-1}
\]
over $F$ at $p\ne\ell$. At primes of good reduction with $p\nmid q\ell$, it is the image of $L_p(A,\chi(p)X)$ under $\iota_\ell$. To compare complex Euler products it suffices to match these factors. The remaining finitely many complex factors may be chosen arbitrarily, provided the factors and their inverses converge absolutely sufficiently far to the right. By \cref{prop:Euler-modifications}, such choices do not affect linear or algebraic independence over $\calF'$ or $\RP$. The same proposition applies to the finitely many factors at which a primitive automorphic twist, when available, may differ from \eqref{eq:twist-definition}.
\end{remark}

\section{Abelian varieties and applications}
\label{sec:applications}

\subsection{Linear independence over \texorpdfstring{$\calF'$}{F'}}

Let $A_1,\ldots,A_N$ be nonzero abelian varieties over $\QQ$, and let $\chi_i$ be Dirichlet characters. Choose a prime $\ell$, a number field $E\subset\CC$ containing their values, and an $\ell$-adic field $F$ (a finite extension of $\QQ_\ell$) with an embedding $\iota_\ell:E\hookrightarrow F$. When we impose a connectedness hypothesis, we choose $\ell$ as in that hypothesis. Write
\[
 L_i(s)=L(s,A_i\otimes\chi_i),\qquad
 W_i=(V_\ell A_i\otimes_{\QQ_\ell}F)\otimes_F\chi_{i,\ell}.
\]
The sets $D_i$ and $D_{ij}$ in \eqref{eq:local-difference-sets} are now formed using these $L_i$.

\Needspace{12\baselineskip}
\begin{theorem}
\label{thm:representation-criterion}
The following conditions are equivalent:
\begin{enumerate}[label=\textup{(\alph*)}]
\item The representations $W_i$ are pairwise non-isomorphic.
\item The family $1,L_1,\ldots,L_N$ is linearly independent over $\calF'$.
\item For every $m\ge0$, the indexed family
\[
 1,\qquad L_i^{(k)}\quad(1\le i\le N,\ 0\le k\le m)
\]
is linearly independent over $\calF'$.
\end{enumerate}
\end{theorem}

\begin{proof}
To prove \textup{(a)}$\Rightarrow$\textup{(c)}, assume that the $W_i$ are pairwise non-isomorphic. For each $i\ne j$, \cref{lem:trace-density} and \eqref{eq:twisted-prime-trace} give a subset of $D_{ij}$ of positive natural density. Every $D_i$ is cofinite, by the degree of the good local polynomial. Apply \cref{cor:density-criteria}\textup{(a)}$\Rightarrow$\textup{(c)}, with $P$ the set of all primes, to obtain \textup{(c)} here. The implication \textup{(c)}$\Rightarrow$\textup{(b)} follows by taking $m=0$.

For \textup{(b)}$\Rightarrow$\textup{(a)}, suppose instead that $W_i\simeq W_j$ for some $i\ne j$. The local factors then agree outside finitely many primes, so $D_{ij}$ is finite. Condition \textup{(a)} of \cref{prop:common} fails for the admissible family $\mathscr T_0$. The equivalence with \textup{(b)} of that proposition gives a nontrivial relation among $1,L_1,\ldots,L_N$, contradicting \textup{(b)} here.
\end{proof}

\begin{proof}[Proof of \cref{thm:main}]
\textup{(1)} Assume that the $A_i$ are pairwise non-isogenous over $\overline\QQ$. If $W_i\simeq W_j$ for some $i\ne j$, restrict this isomorphism to a finite extension $K/\QQ$ on which both characters are trivial. The untwisted Tate modules are then isomorphic after scalar extension to $F$. The equality of traces descends to $\QQ_\ell$, so Brauer--Nesbitt and Faltings's theorem give a $K$-isogeny between $A_i$ and $A_j$, a contradiction. Thus the $W_i$ are pairwise non-isomorphic, and \cref{thm:representation-criterion}\textup{(a)}$\Rightarrow$\textup{(c)} proves \textup{(1)}.

\textup{(2)} Take every $\chi_i$ to be the trivial character modulo $1$, so that $L_i(s)=L(s,A_i)$. By Faltings's theorem, the Tate modules are pairwise non-isomorphic if and only if the varieties are pairwise non-isogenous over $\QQ$. These are precisely the equivalent conditions in \cref{thm:representation-criterion}, so both assertions of \textup{(2)} follow. Explicitly, a $\QQ$-isogeny gives $L(s,A_i)=L(s,A_j)$, including equality of the bad local factors, and hence a nontrivial linear relation.
\end{proof}

\begin{remark}
\label{rem:several-twists}
The same abelian variety may occur with several different twists in \cref{thm:representation-criterion}. Take $A_i=A$ for every $i$ and put $V=V_\ell A\otimes_{\QQ_\ell}F$. The theorem says that $1$, the functions $L(s,A\otimes\chi_i)$, and all their derivatives are linearly independent over $\calF'$ if and only if
\[
 V\not\simeq V\otimes_F(\chi_{j,\ell}\chi_{i,\ell}^{-1})
 \qquad(i\ne j).
\]
In particular, distinct primitive characters suffice when $V$ has no nontrivial finite-order self-twist. If $V$ does have such a self-twist, distinct primitive characters can give isomorphic twisted Tate modules, and independence can fail. Characters inducing the same primitive character also give dependence, since their coefficientwise twists differ at only finitely many Euler factors.
\end{remark}

\subsection{Connected monodromy and completely split primes}

For a finite Galois extension $K/\QQ$, let $\operatorname{Spl}(K)$ be the set of rational primes that split completely in $K$.

\begin{proof}[Proof of \cref{thm:connected}]
We prove \textup{(a)}$\Rightarrow$\textup{(c)}. Assume that the $A_i$ are pairwise non-isogenous over $\overline\QQ$, and fix $P$ with $\ud(P)>0$. Every $D_i$ is cofinite. If $\dim A_i\ne\dim A_j$, the good local polynomials have different degrees, so $D_{ij}$ is also cofinite. If the dimensions agree, Faltings's theorem shows that $V_\ell A_i$ and $V_\ell A_j$ do not become isomorphic over any finite extension. The individual monodromy groups are connected, so \cref{lem:Rajan-density} applies and gives natural density zero for the set of equal twisted traces. Equality of local factors implies equality of traces outside finitely many primes. Hence
\[
 \delta(D_i)=1\quad\text{for every }i,\qquad
 \delta(D_{ij})=1\quad\text{for every }i\ne j.
\]
By \cref{cor:density-criteria}\textup{(a)}$\Rightarrow$\textup{(c)}, the required family is independent over $\RP$. Since $P$ was arbitrary, this proves \textup{(c)}. Taking $m=0$ gives \textup{(c)}$\Rightarrow$\textup{(b)}.

To prove \textup{(b)}$\Rightarrow$\textup{(a)}, assume independence over $\RP$ for every $P$ with $\ud(P)>0$. Suppose that $A_i,A_j$ are geometrically isogenous for some $i\ne j$. Choose a finite Galois extension $K'/\QQ$ over which an isogeny is defined and both characters are trivial. The twisted local factors agree at all but finitely many primes splitting completely in $K'$, so
\[
 D_{ij}\cap\operatorname{Spl}(K')\quad\text{is finite}.
\]
For $P=\operatorname{Spl}(K')$ we have $\delta(P)=1/[K':\QQ]>0$ and $\ud(P\cap D_{ij})=0$. Thus condition \textup{(a)} of \cref{cor:density-criteria} fails. Its equivalence with \textup{(b)} gives dependence of $1,L_1,\ldots,L_N$ over this $\RP$, a contradiction. This proves \textup{(a)}.
\end{proof}

\begin{corollary}
\label{cor:splitting}
Let $A_1,\ldots,A_N$ be nonzero abelian varieties over $\QQ$ that are pairwise non-isogenous over $\overline\QQ$, and let $\chi_i$ be arbitrary Dirichlet characters. There is a finite Galois extension $K/\QQ$ such that, for every set $P$ of primes with
\[
 \ud\bigl(P\cap\operatorname{Spl}(K)\bigr)>0,
\]
the family consisting of $1$ and all derivatives of $L(s,A_i\otimes\chi_i)$ is linearly independent over $\RP$. One may take $K$ to contain the connectedness field of $\bigoplus_i V_\ell A_i$ and the fixed fields of all $\chi_{i,\ell}$.
\end{corollary}

\begin{proof}
Let $G$ be the monodromy group of $\bigoplus_i V_\ell A_i$ and choose $K$ as stated. The image of $G_K$ is Zariski dense in $G^0$, and the characters are trivial on $G_K$. For $i\ne j$, let $t_{ij}$ be the difference of the two trace functions on $G$. Its restriction to $G^0$ is nonzero: otherwise Brauer--Nesbitt, semisimplicity over $K$, and Faltings would give a $K$-isogeny.

Apply \cref{thm:algebraic-chebotarev} to
\[
 X_{ij}=G^0\cap\{t_{ij}=0\}\subseteq G.
\]
This is a conjugation-stable closed subvariety containing no connected component of $G$. For all but finitely many $p\in\operatorname{Spl}(K)$, the characters have value $1$ and the local polynomials are given by \eqref{eq:prime-trace}. At such a prime, equality of the local factors implies $t_{ij}(\Frob_p)=0$. It follows that
\[
 \delta\bigl(\operatorname{Spl}(K)\setminus D_{ij}\bigr)=0.
\]
Consequently $\ud(P\cap D_{ij})>0$ under the stated hypothesis. The $D_i$ are cofinite, so \cref{cor:density-criteria} applies.
\end{proof}

\begin{corollary}
\label{cor:partial-products}
Assume all the hypotheses of \cref{thm:connected}, including connectedness of the individual monodromy groups, and its condition \textup{(a)}. Let $P$ have positive upper density and put
\[
 L_{i,P}(s)=\prod_{p\in P}L_p(A_i,\chi_i(p)p^{-s}).
\]
Then $1$ and all derivatives of the $L_{i,P}$ are linearly independent over $\RP$. The same conclusion holds under the hypotheses on $A_i$ and $P$ in \cref{cor:splitting}.
\end{corollary}

\begin{proof}
Let $E_{i,P^c}$ be the product of the Euler factors of $L_i$ at primes outside $P$. It is a unit of $\RP$, and $L_i=E_{i,P^c}L_{i,P}$. Apply \cref{prop:Euler-modifications}.
\end{proof}

\Needspace{11\baselineskip}
\begin{corollary}
\label{cor:necessary-density}
Assume all the hypotheses of \cref{thm:connected}, including connectedness of the individual monodromy groups, and its condition \textup{(a)}. Suppose that a nontrivial relation
\[
 C_0+\sum_{i=1}^N\sum_{k=0}^m C_{ik}L_i^{(k)}=0,
 \qquad C_0,C_{ik}\in\Dabs,
\]
holds, and let $T$ be the union of the prime supports of its coefficients. Then $T$ contains all but finitely many primes in at least one of the sets $D_i$ or $D_{ij}$. In particular, $T$ has natural density one.
\end{corollary}

\begin{proof}
The first assertion is \cref{cor:relation-support}. Every $D_i$ contains all but finitely many primes, and the proof of \cref{thm:connected} shows that every $D_{ij}$ has natural density one. Hence $T$ has natural density one.
\end{proof}

\begin{remark}
\label{rem:quadratic-twists}
The quantifier over all $P$ in \cref{thm:connected} is essential. Let $E/\QQ$ be an elliptic curve for which the monodromy group of $V_\ell E$ is $\GL_2$. Let $E^\eta$ be its twist by a nontrivial quadratic character $\eta$. Both individual monodromy groups are $\GL_2$. They are geometrically isogenous but not $\QQ$-isogenous: an intertwiner would commute with the Zariski-dense image of $\ker\eta$, hence be scalar, which is incompatible with the nontrivial twist. The functions $1,L(s,E),L(s,E^\eta)$ are linearly independent over $\calF'$ by \cref{thm:representation-criterion}. For
\[
 P=\{p:\eta(p)=1\},
\]
of density $1/2$, the local factors agree at all but finitely many primes of $P$; equivalently, $D_{12}\cap P$ is finite. Thus the equivalence of \textup{(a)} and \textup{(b)} in \cref{cor:density-criteria} gives dependence of $1,L(s,E),L(s,E^\eta)$ over $\RP$, illustrating \textup{(b)}$\Rightarrow$\textup{(a)} of \cref{thm:connected}.

Even the assumption that the varieties are pairwise non-isogenous over $\overline\QQ$ does not suffice for arbitrary $P$ without a monodromy assumption. Take elliptic curves $E_1,E_2/\QQ$ with CM fields $\QQ(i)$ and $\QQ(\sqrt{-2})$. They are geometrically non-isogenous. Every good prime $p\equiv7\pmod8$ is inert in both CM fields, so both traces vanish and both Euler factors are $(1+pX^2)^{-1}$. Taking $P$ to be this progression, of density $1/4$, again gives dependence over $\RP$. This does not contradict \cref{cor:splitting}, which restricts the sets $P$.

The quadratic twist $E^\eta$ above also shows that the monodromy group of a direct sum can become disconnected when one summand is twisted. For $V=V_\ell E$, the monodromy groups of $V\oplus V$ and $V\oplus(V\otimes\eta)$ are
\[
 \{(g,g):g\in\GL_2\},\qquad
 \{(g,g):g\in\GL_2\}\cup\{(g,-g):g\in\GL_2\},
\]
respectively. The first is connected and the second has two components.
\end{remark}

\subsection{Algebraic independence}

\begin{proof}[Proof of \cref{thm:algebraic}]
We prove \textup{(a)}$\Rightarrow$\textup{(b)}$\Rightarrow$\textup{(c)}$\Rightarrow$\textup{(d)}$\Rightarrow$\textup{(a)}. The implication \textup{(a)}$\Rightarrow$\textup{(b)} follows from $\CC\subseteq\calF'$.

To prove \textup{(b)}$\Rightarrow$\textup{(c)}, suppose that \textup{(c)} fails. Multiply a nonzero rational relation among the multiplicity vectors by a common denominator and write the resulting integer vector as $\alpha-\beta$, with $\alpha,\beta\in\ZZ_{\ge0}^r$ distinct. Then
\[
 \sum_i\alpha_i n_{ij}=\sum_i\beta_i n_{ij}\quad(1\le j\le t),
 \qquad
 \prod_i A_i^{\alpha_i}\sim_\QQ\prod_i A_i^{\beta_i}.
\]
The $L$-functions of these products are equal, so
\[
 \prod_i L(s,A_i)^{\alpha_i}
 -\prod_i L(s,A_i)^{\beta_i}=0.
\]
This gives a nontrivial polynomial relation over $\CC$, contradicting \textup{(b)}.

For \textup{(c)}$\Rightarrow$\textup{(d)}, Faltings's theorem gives
\[
 \operatorname{Hom}_{G_\QQ}(V_\ell B_j,V_\ell B_k)=0\quad(j\ne k).
\]
Since these Tate modules are semisimple, distinct $B_j$ have disjoint sets of irreducible constituents. Their classes are therefore linearly independent in $K_0(\QQ_\ell)\otimes_\ZZ\QQ$. The identities
\[
 [V_\ell A_i]=\sum_{j=1}^t n_{ij}[V_\ell B_j]
 \quad(1\le i\le r)
\]
show that \textup{(c)} implies \textup{(d)}.

To prove \textup{(d)}$\Rightarrow$\textup{(a)}, assume \textup{(d)}. For distinct exponent vectors $\alpha,\beta\in\ZZ_{\ge0}^r$, we have
\[
 \sum_i\alpha_i[V_\ell A_i]\ne\sum_i\beta_i[V_\ell A_i].
\]
Thus the representations
\[
 V_\alpha=\bigoplus_{i=1}^r(V_\ell A_i)^{\oplus\alpha_i}
 \quad(\alpha\in\ZZ_{\ge0}^r\setminus\{0\})
\]
are nonzero and pairwise non-isomorphic. They are the Tate modules of the products $A_\alpha=\prod_i A_i^{\alpha_i}$, whose $L$-functions are
\[
 M_\alpha(s)=L(s,A_\alpha)=\prod_{i=1}^r L(s,A_i)^{\alpha_i}.
\]
The zero exponent vector gives $M_0=1=L(s,0)$. For each finite collection of nonzero exponent vectors, \cref{thm:representation-criterion}, with $F=\QQ_\ell$ and trivial characters modulo $1$, gives linear independence of the corresponding $M_\alpha$, together with $1$, over $\calF'$. This excludes every polynomial relation and proves \textup{(a)}.

Finally, assume \textup{(d)} and suppose that the monodromy group $G$ of $\bigoplus_i V_\ell A_i$ is connected. View $G$ as a subgroup of $\prod_i\GL(V_\ell A_i)$. Fix a finite nonempty set
$\Lambda\subset\ZZ_{\ge0}^r\setminus\{0\}$ of exponent vectors. The monodromy group $G_\Lambda$ of $\bigoplus_{\alpha\in\Lambda}V_\alpha$ is the image of $G$ under
\[
 (g_1,\ldots,g_r)\longmapsto
 \left(\bigoplus_{i=1}^r g_i^{\oplus\alpha_i}\right)_{\alpha\in\Lambda}
 \in\prod_{\alpha\in\Lambda}\GL(V_\alpha).
\]
Thus $G_\Lambda$ is connected.

Choose a finite set $S$ containing $\ell$ and all primes of bad reduction of the $A_i$. The $V_\alpha$ are pairwise non-isomorphic by \textup{(d)}, so \cref{lem:joint-density}, with $F=\QQ_\ell$, gives
\[
 \delta\left(\left\{p\notin S:
 \Tr(\Frob_p\mid V_\alpha)=\Tr(\Frob_p\mid V_\beta)
 \right\}\right)=0
 \qquad(\alpha\ne\beta\text{ in }\Lambda).
\]

Let $D_\alpha$ and $D_{\alpha\beta}$ be the sets in \eqref{eq:local-difference-sets} for the Euler products $M_\alpha$. Equal local factors have equal prime coefficients, hence equal Frobenius traces outside $S$. Therefore $\delta(D_{\alpha\beta})=1$. Also, $D_\alpha$ contains every prime outside $S$, since the local polynomial has degree $\dim V_\alpha>0$. For every $P$ with $\ud(P)>0$, we obtain
\[
 \begin{aligned}
 \ud(P\cap D_\alpha)&=\ud(P)>0
 &&(\alpha\in\Lambda),\\
 \ud(P\cap D_{\alpha\beta})&=\ud(P)>0
 &&(\alpha\ne\beta\text{ in }\Lambda).
 \end{aligned}
\]
By \cref{cor:density-criteria}, the family $1,(M_\alpha)_{\alpha\in\Lambda}$ is linearly independent over $\RP$. Since $\Lambda$ was arbitrary, this proves algebraic independence over every such $\RP$. Conversely, algebraic independence over any $\RP$ implies (a), since $\calF'\subseteq\RP$.
\end{proof}

\begin{remark}
\label{rem:trace-reduction}
The criterion in \cref{thm:algebraic} can also be expressed through Frobenius traces modulo $p$. Put
\[
 \mathcal A=
 \left(\prod_p\mathbb F_p\right)\Big/
 \left(\bigoplus_p\mathbb F_p\right),
\]
where $p$ runs over the rational primes. The map $\ZZ\to\mathcal A$ given by $n\mapsto(n\bmod p)_p$ is injective. Every nonzero integer has an invertible image, since only finitely many primes divide it. This embedding therefore extends uniquely to $\QQ\hookrightarrow\mathcal A$. This is also the $\QQ$-algebra in which Kaneko--Zagier's finite multiple zeta values take their values; see \cite[\S7]{Kaneko}. For an abelian variety $A/\QQ$, define
\[
 \mathbf a(A)=\bigl(a_A(p)\bmod p\bigr)_p\in\mathcal A,
\]
with arbitrary components at the finitely many primes of bad reduction.

Let $K_0^{\mathrm{AV}}\subseteq K_0(\QQ_\ell)$ be the subgroup generated by the classes $[V_\ell A]$ as $A$ varies over abelian varieties over $\QQ$. For integers $m_i$, a relation among these classes gives the corresponding relation among Frobenius traces at all sufficiently large primes:
\[
 \sum_i m_i[V_\ell A_i]=0
 \quad\Longrightarrow\quad
 \sum_i m_i a_{A_i}(p)=0
 \quad\text{for all sufficiently large }p.
\]
Reducing modulo $p$ therefore defines a $\QQ$-linear map
\[
 \Phi:K_0^{\mathrm{AV}}\otimes_\ZZ\QQ\longrightarrow\mathcal A,
 \qquad
 [V_\ell A]\longmapsto\mathbf a(A).
\]
This map is injective. Let $x\in\ker\Phi$. Choose a positive integer $d$ and integers $m_i$ such that
\[
 d x=\sum_i m_i[V_\ell A_i],
 \qquad
 \sum_i m_i\mathbf a(A_i)=0.
\]
For all sufficiently large primes of common good reduction,
\[
 p\mid t_p,\qquad
 t_p=\sum_i m_i a_{A_i}(p),\qquad
 |t_p|\le
 2\left(\sum_i|m_i|\dim A_i\right)\sqrt p<p.
\]
The bound is the Weil bound for Frobenius traces. Thus $t_p=0$ at all but finitely many primes. Applying \cref{lem:trace-density} to the positive and negative sums of Tate modules shows that they are isomorphic. Hence
$\sum_i m_i[V_\ell A_i]=0$ in $K_0(\QQ_\ell)$. Thus $d x=0$, and hence $x=0$, proving injectivity.

Consequently, conditions \textup{(a)}--\textup{(d)} of \cref{thm:algebraic} are also equivalent to $\QQ$-linear independence of
\[
 \mathbf a(A_1),\ldots,\mathbf a(A_r)\quad\text{in }\mathcal A.
\]
In particular, these elements are linearly independent when the $A_i$ are $\QQ$-simple and pairwise non-isogenous over $\QQ$.
\end{remark}

\Needspace{16\baselineskip}
\begin{proposition}
\label{prop:twisted-algebraic}
Let $A_1,\ldots,A_N$ be nonzero abelian varieties over $\QQ$, and let $\chi_1,\ldots,\chi_N$ be Dirichlet characters. Choose a prime $\ell$, a number field $E\subset\CC$ containing their values, and an embedding $\iota_\ell:E\hookrightarrow F$ into an $\ell$-adic field $F$. Let $\chi_{i,\ell}:G_\QQ\to F^\times$ be the finite-order Galois character associated, via $\iota_\ell$, with the primitive character inducing $\chi_i$. Define
\[
 L_i(s)=L(s,A_i\otimes\chi_i),\qquad
 W_i=(V_\ell A_i\otimes_{\QQ_\ell}F)\otimes_F\chi_{i,\ell}.
\]
Let $K_0(F)$ be the Grothendieck group of continuous semisimple representations of $G_\QQ$ on finite-dimensional $F$-vector spaces. The following conditions are equivalent:
\begin{enumerate}[label=\textup{(\alph*)}]
\item $L_1,\ldots,L_N$ are algebraically independent over $\calF'$.
\item The classes $[W_1],\ldots,[W_N]$ are linearly independent in $K_0(F)\otimes_\ZZ\QQ$.
\end{enumerate}
If the monodromy group of $\bigoplus_i W_i$ is connected, these conditions are also equivalent to algebraic independence of $L_1,\ldots,L_N$ over $\RP$ for every set $P$ of primes with $\ud(P)>0$.
\end{proposition}

\begin{proof}
To prove \textup{(b)}$\Rightarrow$\textup{(a)}, assume \textup{(b)}. For distinct nonzero exponent vectors $\alpha$, the semisimple representations
\[
 W_\alpha=\bigoplus_i W_i^{\oplus\alpha_i}
\]
are nonzero and pairwise non-isomorphic. Put $M_\alpha=\prod_i L_i^{\alpha_i}$. Away from finitely many primes, the Euler factors attached to $W_\alpha$ agree with those of $M_\alpha$. For every pair, \cref{lem:trace-density} gives a set of primes of positive density on which the traces, and hence the local factors, differ. Each $M_\alpha$ has a local polynomial of positive degree outside finitely many primes. Thus \cref{cor:density-criteria}, with $P$ the set of all primes, gives linear independence of every finite collection of the $M_\alpha$, together with $1$, over $\calF'$. This proves \textup{(a)}.

To prove \textup{(a)}$\Rightarrow$\textup{(b)}, suppose that \textup{(b)} fails. Multiply a nonzero rational relation among the classes $[W_i]$ by a common denominator. Separating positive and negative coefficients gives $W_\alpha\simeq W_\beta$ for distinct exponent vectors $\alpha,\beta\in\ZZ_{\ge0}^N$, since the representations are semisimple. The monomials then have equal local factors outside a finite set $S$. Let $E_{\alpha,S}$ be the product of the Euler factors of $M_\alpha$ at primes in $S$, and define $E_{\beta,S}$ in the same way. Then
\[
 M_\alpha=U M_\beta,\qquad
 U=\frac{E_{\alpha,S}}{E_{\beta,S}}\in\calF'^\times.
\]
Since $\alpha\ne\beta$, the polynomial $X^\alpha-U X^\beta$ is nonzero. The displayed relation therefore contradicts \textup{(a)} and proves \textup{(a)}$\Rightarrow$\textup{(b)}.

Finally, assume \textup{(b)} and connectedness of the monodromy group of $\bigoplus_i W_i$. For each finite collection of exponent vectors, the monodromy group of the corresponding direct sum of the $W_\alpha$ is an image of this connected group, hence is connected. By \cref{lem:joint-density}, the set of primes with equal traces for any pair has natural density zero. Applying \cref{cor:density-criteria} as above gives algebraic independence over every $\RP$ with $\ud(P)>0$. The converse follows from $\calF'\subseteq\RP$.
\end{proof}

The untwisted equivalence with algebraic independence over $\CC$ in \cref{thm:algebraic} uses exact equality of the $L$-functions of $\QQ$-isogenous abelian varieties. In the coefficientwise twisted case, the ratio $U$ of the two finite products of Euler factors need not be constant. The relation proved above therefore has coefficients in $\calF'$, and it need not give a relation over $\CC$.

\subsection{Curves and Jacobians}

Let $C/\QQ$ be a smooth projective geometrically connected curve of positive genus. The Abel--Jacobi map gives
\[
 H^1_{\mathrm{\acute et}}(C_{\overline\QQ},\QQ_\ell)
 \simeq H^1_{\mathrm{\acute et}}(J(C)_{\overline\QQ},\QQ_\ell)
 \simeq(V_\ell J(C))^\vee;
\]
see \cite[Corollary~9.6 and Remark~11.5]{MilneJac}. Changing the geometric base point changes the map by a translation, which acts trivially on $H^1$. Hence the induced isomorphism is Galois-equivariant. For the $H^1$-Hasse--Weil $L$-function of $C$ and its coefficientwise Dirichlet twist, we write
\[
 L(s,C)=L(s,J(C)),\qquad
 L(s,C\otimes\chi)=L(s,J(C)\otimes\chi).
\]

\begin{corollary}
\label{cor:curves}
Let $C_1,\ldots,C_N$ be smooth projective geometrically connected curves over $\QQ$ of positive genus, put $J_i=J(C_i)$, and fix Dirichlet characters $\chi_1,\ldots,\chi_N$.
\begin{enumerate}[label=\textup{(\arabic*)}]
\item The $J_i$ are pairwise non-isogenous over $\QQ$ if and only if, for every $m\ge0$, the indexed family
\[
 1,\qquad L^{(k)}(s,C_i)
 \quad(1\le i\le N,\ 0\le k\le m)
\]
is linearly independent over $\calF'$.
\item If the $J_i$ are pairwise non-isogenous over $\overline\QQ$, then for every $m\ge0$ the indexed family
\begin{equation}
\label{eq:curve-twisted-family}
 1,\qquad L^{(k)}(s,C_i\otimes\chi_i)
 \quad(1\le i\le N,\ 0\le k\le m)
\end{equation}
is linearly independent over $\calF'$.
\item Suppose that, for some prime $\ell$, the monodromy group of each $V_\ell J_i$ is connected. Then
the $J_i$ are pairwise non-isogenous over $\overline\QQ$ if and only if, for every set $P$ of primes with $\ud(P)>0$ and every $m\ge0$, the family \eqref{eq:curve-twisted-family} is linearly independent over $\RP$.
\item If the $J_i$ are pairwise non-isogenous over $\overline\QQ$, there is a finite Galois extension $K/\QQ$ such that, for every set $P$ of primes with
\[
 \ud\bigl(P\cap\operatorname{Spl}(K)\bigr)>0
\]
and every $m\ge0$, the family \eqref{eq:curve-twisted-family} is linearly independent over $\RP$. No connectedness assumption is needed in this assertion.
\end{enumerate}
\end{corollary}

\begin{proof}
Apply \cref{thm:main}\textup{(2)} to $A_i=J_i$ for \textup{(1)}, and \cref{thm:main}\textup{(1)} for \textup{(2)}. Conditions \textup{(3)} and \textup{(4)} follow from \cref{thm:connected} and \cref{cor:splitting}, respectively, using the stated identities of $L$-functions.
\end{proof}

\begin{remark}
\label{rem:curve-hypotheses}
The following sufficient conditions can be used in \cref{cor:curves}.
\begin{enumerate}[label=\textup{(\arabic*)}]
\item Different genera imply that the Jacobians are not isogenous over $\overline\QQ$. For each pair of curves of the same genus, a single prime $p$ of common good reduction with
\[
 P_p(J_i,X)\ne P_p(J_j,X)
\]
rules out an isogeny over $\QQ$. In particular, different numbers of $\mathbb F_p$-points on the two reductions suffice. These tests do not, in general, rule out an isogeny over $\overline\QQ$.

One can also rule out an isogeny over $\overline\QQ$ by a local test. At a prime $p$ of common good reduction, suppose that $J_i$ has a Frobenius eigenvalue $\alpha$ such that $\alpha/\beta$ is not a root of unity for any Frobenius eigenvalue $\beta$ of $J_j$. Then $J_i$ and $J_j$ are not isogenous over $\overline\QQ$. Indeed, an isogeny over a number field would force the multisets of $d$-th powers of their Frobenius eigenvalues at $p$ to agree for some $d\ge1$, contradicting this condition.

\item Suppose that $C_i$ is the smooth projective model of $y^2=f_i(x)$, where $f_i\in\QQ[x]$ is irreducible of degree $n_i\ge5$ and its Galois group is the symmetric group $S_{n_i}$ or the alternating group $A_{n_i}$. Then
\[
 \operatorname{End}_{\overline\QQ}(J_i)=\ZZ
\]
by \cite[Theorem~2.1]{ZarhinEnd}. If the splitting fields of the $f_i$ are pairwise linearly disjoint over $\QQ$, the $J_i$ are pairwise non-isogenous over $\overline\QQ$ by \cite[Theorem~1.2]{ZarhinHom}. This verifies the hypothesis that the Jacobians are pairwise non-isogenous over $\overline\QQ$ in \cref{cor:curves}\textup{(2)--(4)}.

\item In genus one, a non-CM Jacobian has monodromy group $\GL_{2,\QQ_\ell}$ for every $\ell$, by Serre's open-image theorem \cite{SerreOpenImage}. In higher genus, let $J/\QQ$ be a Jacobian of dimension $g\ge2$ with $\operatorname{End}_{\overline\QQ}(J)=\ZZ$. Suppose that $J$ has semistable reduction of toric rank one at some prime. Let $H_\ell\subseteq\operatorname{GSp}_{2g}(\ZZ_\ell)$ be the Galois image on $T_\ell J$. Hall's theorem \cite[Theorem~1]{Hall} shows that its reduction modulo $\ell$ is $\operatorname{GSp}_{2g}(\mathbb F_\ell)$ for all sufficiently large $\ell$. For $\ell\ge5$, the closure of the commutator subgroup of $H_\ell$ surjects onto $\operatorname{Sp}_{2g}(\mathbb F_\ell)$, so it equals $\operatorname{Sp}_{2g}(\ZZ_\ell)$ by \cite[Theorem~1]{LiftingSymplectic}. Since the similitude character is cyclotomic and has infinite image, the Zariski closure of $H_\ell$ is $\operatorname{GSp}_{2g,\QQ_\ell}$, which is connected.

For the hyperelliptic curves in \textup{(2)}, this reduction condition holds if $f_i\in\ZZ[x]$ is monic and, for some odd prime $p_i$, its reduction has the form
\[
 \overline f_i(X)=(X-a_i)^2h_i(X),\qquad
 h_i(a_i)\ne0,\qquad h_i\text{ squarefree},
\]
with $a_i\in\mathbb F_{p_i}$ and $h_i\in\mathbb F_{p_i}[X]$; see \cite[Introduction]{Hall}. Combining this condition for every $i$ with the disjointness condition in \textup{(2)} gives Jacobians that are pairwise non-isogenous over $\overline\QQ$ and whose individual monodromy groups are connected, as required in \cref{cor:curves}\textup{(3)}. A common sufficiently large prime $\ell$ works for the whole family.
\end{enumerate}
\end{remark}

\bibliographystyle{amsplain}
\bibliography{hasse_weil_references}
\end{document}